\documentclass[preprint]{elsarticle}

\usepackage{graphicx}
\usepackage{amssymb,amscd}
\usepackage{amsmath,amssymb}
\usepackage{bbm}
\usepackage{epsfig}
\usepackage[left=2cm,top=2cm,right=2cm,bottom=2cm]{geometry}
\usepackage{enumitem}
\setitemize{listparindent=\parindent}

\begin{document}

\title{Adjoint method for a tumor invasion PDE-constrained optimization problem in 2D using Adaptive Finite Element Method.}

\author[famaf,ciem]{A.A.I.~Quiroga \corref{cor1}}
\ead{aiquiroga@famaf.unc.edu.ar}

\author[famaf,ciem]{D.R.~Fern\'andez}
\ead{dfernandez@famaf.unc.edu.ar}

\author[famaf,ciem]{G.A.~Torres}
\ead{torres@famaf.unc.edu.ar}

\author[famaf,ciem]{C.V.~Turner} 
\ead{turner@famaf.unc.edu.ar} 

\cortext[cor1]{Corresponding author}

\address[famaf]{Facultad de Matem\'atica, Astronom\'ia y F\'isica, Medina Allende s/n, 5000 C\'ordoba, Argentina}

\address[ciem]{Centro de Investigaciones y Estudios en Matem\'atica - CONICET, Medina Allende s/n, 5000 C\'ordoba, Argentina}

\begin{abstract}
In this paper we present a method for estimating unknown parameter that appear in a two dimensional non-linear reaction-diffusion model of cancer invasion. This model considers that tumor-induced alteration of micro-environmental pH provides a mechanism for cancer invasion. A coupled system reaction-diffusion describing this model is given by three partial differential equations for the 2D non-dimensional spatial distribution and temporal evolution of the density of normal tissue, the neoplastic tissue growth and the excess concentration of H$^+$ ions. Each of the model parameters has a corresponding biological interpretation, for instance, the growth rate of neoplastic tissue, the diffusion coefficient, the re-absorption rate and the destructive influence of H$^+$ ions in the healthy tissue.
  
After solving the direct problem, we propose a model for the estimation of parameters by fitting the numerical solution with real data, obtained via \emph{in vitro} experiments and fluorescence ratio imaging microscopy. We define an appropriate functional to compare both the real data and the numerical solution using the adjoint method for the minimization of this functional.

We apply a splitting strategy joint with Adaptive Finite Element Method (AFEM) to solve the direct problem and the adjoint problem. The minimization problem (the inverse problem) is solved by using a trust-region-reflective method including the computation of the derivative of the functional.
\end{abstract}


\maketitle

\section{Introduction}
Cancer is one of the diseases causing the most deaths in the world, despite the best efforts of medicine. Human and financial resources are devoted for cancer research, and on several occasions these efforts are successful \cite{Adam1, Adambellomo, BeChaDe09, BellomoLiMaini, byrne2010dissecting, bellomo2000modelling}.

Some comments on the importance of mathematical modeling in cancer can be found in the literature. In the work \cite{BellomoLiMaini} the authors say ``Cancer modelling has, over the years, grown immensely as one of the challenging topics involving applied mathematicians working with researchers active in the biological sciences. The motivation is not only scientific as in the industrial nations cancer has now moved from seventh to second place in the league table of fatal diseases, being surpassed only by cardiovascular diseases.''

We use the analysis proposed by Gatenby in \cite{gatenby1996reaction}, which supports the acid-mediated invasion hypothesis. Therefore, it can be represented mathematically as a reaction-diffusion system which describes the spatial and temporal evolution of the tumor tissue, normal tissue, and excess concentration of H$^+$.

The model simulates a pH gradient extending from the tumor-host interface. The effect of biological parameters that control this transition is supported by experimental and clinical observations \cite{martin1994noninvasive}.

Some authors \cite{gatenby1996reaction} model tumor invasion in order to find an underlying mechanism by which primary and metastatic cancers invade and destroy normal tissues. They do not attempt to model the genetic changes that lead to the transformation and seek to understand the causes of these changes. Likewise, they do not attempt to model the large-scale morphological aspects of tumor necrosis such as central necrosis. Instead, they concentrate on the interactions of microscopic scale populations that occur at the tumor-host interface, arguing that these processes influence the clinically significant manifestations of invasive cancer.

Moreover, in \cite{gatenby1996reaction}, the authors suppose that transformation-induced reversion of neoplastic tissue to primitive glycolytic metabolic pathways, with resultant increased acid production and the diffusion of that acid into surrounding healthy tissue, creates a peritumoral micro-environment in which the tumor cells survive and proliferate, while normal cells may not remain viable. The following temporal sequence would derive: (a) a high concentration of H$^+$ ions in tumors will diffuse chemically as a gradient to adjacent normal tissue, exposing these normal cells to an interstitial pH like in the tumor, (b) normal cells immediately adjacent to the edge of the tumor are unable to survive in chronically this acid environment, and (c) progressive loss of normal cell layers in the tumor-host interface facilitates tumor invasion. Key elements of this mechanism of tumor invasion include low pH due to primitive metabolism and reduced viability of normal tissue in an acidic environment.

This model depends only on a small number of cellular and sub-cellular parameters. The analysis of the equations shows that the model simulates a crossover from a benign tumor to a malignant invasive tumor when some combination of parameters turn over some threshold value.

\begin{figure}
\begin{center}
\includegraphics[scale=0.3]{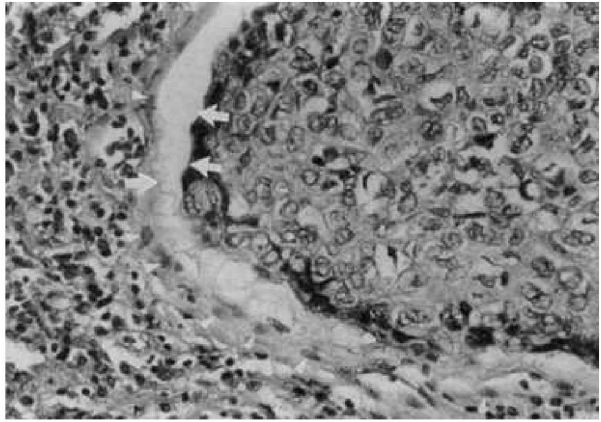}
\setlength{\belowcaptionskip}{2pt}
\caption{A micrographs of the tumor-host interface from human squamous cell carcinomas of the head and neck \cite{gatenby1996reaction}.}
\label{fig:gap}
\end{center}
\end{figure}

The structure and dynamics of the tumor-host interface in invasive cancers are controlled by the same parameters which generate a transformation from a benign tumor into malignant tumor. A hypo-cellular interstitial space, as we can see in Figure \ref{fig:gap} \cite[Figure 4a]{gatenby1996reaction}, occurs in some cancers.

In this paper we propose a framework via a PDE-constrained optimization problem that can be solved with the splitting method, exploiting the fact that this procedure is easily parallelizable. We follow the PDE-based model by Gatenby \cite{gatenby1996reaction} in a two-dimensional tissue. We estimate one of the model parameters (the destructive influence of H$^+$ ions in the healthy tissue) using an inverse problem. It is possible to get data about the concentration of hydrogen ions \cite{martin1994noninvasive} via fluorescence ratio imaging microscopy. In this approach, tumor invasion is modeled via a coupled nonlinear system of partial differential equations, which makes the numerical solution procedure quite challenging. These equations are solved using an Adaptive Finite Element Method (AFEM). 

This kind of problem constitutes a particular application of the so-called inverse problems, which are being increasingly used in a broad number of fields in applied sciences. For instance, problems referred to structured population dynamics \cite{PeZu07}, computerized tomography and image reconstruction in medical imaging \cite{van2011source,ZuMa03}, and more specifically tumor growth \cite{AgBaTu,hogea,knopoff}, among many others.

We solve a minimization problem using a gradient-based method considering the adjoint method in order to find the derivative of an objective functional. In this way, we would obtain the best parameter that fits patient-specific data.

This work follows the ideas \cite{AQFETOTU-1d} where the space variable was in a one dimensional space. The extension of the model  to two dimensional space allows us to approach the results to more realistic biological hypotheses.

The contents of this paper, which is organized into 5 sections, are as follows: Section 2 consists in some preliminaries about the model, the definition of variational form of the direct and adjoint problems, and the minimization problem. Section 3 deals with suitable numerical algorithms to solve the direct and adjoint problems. In particular, we use the splitting method and the Adaptive Finite Element Method with a computation of \emph{a posteriori} error. In Section 4 we show numerical simulations of the retrieved parameter and the need of a parallel scheme. Section 5 presents the conclusions and some future work related to the contents of this paper.

Some words about our notation.  We use $\langle\cdot,\cdot\rangle$ to denote the $L^2$ inner product (the space is always clear from the context) and we consider the sum of inner products for a Cartesian product of spaces.  For a function $F:V\times U_{ad}\rightarrow \mathcal{Z}$ such that $(u,\delta_1)\mapsto F(u,\delta_1)$, we denote by $F'(u,\delta_1)$ the full Fr\'echet-derivative and by $\frac{\partial F}{\partial u}(u,\delta_1)$ and $\frac{\partial F}{\partial \delta_1}(u,\delta_1)$ the partial Fr\'echet-derivatives of $F$ at $(u,\delta_1)$. For a linear operator $T : V \rightarrow \mathcal{Z}$ we denote $T^* : \mathcal{Z}^* \rightarrow V^*$ the adjoint operator of $T$. If $T$ is invertible, we call $T^{-*}$ the inverse of the adjoint operator $T^*$.

\section{Some preliminaries about the model}
A mathematical model of the tumor-host interface based on the acid mediation hypothesis of tumor invasion due to \cite{gatenby1996reaction} is given by the following system of partial differential equations (PDEs):
\begin{eqnarray*}
\frac{ \partial N_1 }{ \partial t } &=& r_1 N_1 \left( 1 - \frac{N_1 }{K_1} \right) - d_1 L  N_1, \\ 
\frac{ \partial N_2 }{ \partial t } &=& r_2N_2 \left( 1 - \frac{N_2 }{K_2} \right)+ \nabla \cdot \left( D_{N_2} \left( 1 - \frac{N_1}{K_1} \right) \nabla N_2 \right), \\ 
\frac{ \partial L}{ \partial t } &=& r_3N_2 - d_3 L  + D_{N_3} \Delta L, 
\end{eqnarray*}
where the variables are in $\Omega \times [0,T]$. These equations determine the spatial distribution and temporal evolution of three fields: $N_1$, the density of normal tissue; $N_2$, the density of neoplastic tissue; and $L$, the excess concentration of H$^+$ ions. The units of $N_1$ and $N_2$ are cells/cm$^3$ and $L$ is expressed as a molarity (M). The space $x$ and time $t$ are given in cm and seconds, respectively. 

The biological meaning of each equation can be seen in \cite{gatenby1996reaction}, and the non-dimensional mathematical model is:
\begin{eqnarray}
\frac{ \partial u_1 }{ \partial t }  &=& u_1( 1 - u_1 ) - \delta_1u_1u_3, \nonumber \\
\frac{ \partial u_2 }{ \partial t }  &=& \rho_2u_2( 1 - u_2 ) + \nabla \cdot \left( D_2( 1 - u_1) \nabla u_2  \right), \label{nondim} \\
\frac{ \partial u_3 }{ \partial t }  &=& \delta_3( u_2 - u_3 )  + \Delta u_3, \nonumber
\end{eqnarray} 
where the four dimensionless quantities which parameterize the model are given by: $\delta_1 = d_1 r_3 K_2/(d_3 r_1)$, $\rho_2 = r_2/r_1$, $D_2 = D_{N_2}/D_{N_3}$ and $\delta_3 = d_3/r_1$.

The interaction parameters between different cells (healthy and tumor) and concentration of H$^+$ are difficult to measure experimentally. This is the reason for which we propose to estimate them, so we will focus on $\delta_1$ in this work.

The initial and boundary conditions considered for the non-dimensional system are:
\[
\begin{array}{ccc}
u_1(x,0) = u_1^0(x), & u_2(x,0) = u_2^0(x), & u_3(x,0) = u_3^0(x), \quad \forall \, x \in \Omega, \\[2mm]
\displaystyle \frac{ \partial u_1 }{ \partial n }(x,t) = 0, & \displaystyle \frac{ \partial u_2 }{ \partial n }(x,t) = 0, & \displaystyle \frac{ \partial u_3 }{ \partial n }(x,t) = 0, \quad \forall \, x \in \partial \Omega. 
\end{array}
\]

From now on, equations (\ref{nondim}) with the initial and boundary conditions will be referred to as the direct problem.

\subsection{Variational form for the direct problem}
Using the variational techniques for obtaining the weak solution of the direct problem \cite{ladyzhenskaia1988linear,kinderlehrer1987introduction,evans1998partial}, it can be written as $E(u,\delta_1)=0$ where $E : V \times U_{ad} \to V^* $ such that
\begin{eqnarray}
\langle E(u,\delta_1),\lambda\rangle & = & \int_0^T \int_{\Omega} \left( \frac{ \partial u_1 }{ \partial t } \lambda_1 - u_1(1-u_1)\lambda_1 + \delta_1u_1u_3\lambda_1 \right) dxdt + \nonumber \\
&  & \int_0^T \int_{\Omega} \left( \frac{ \partial  u_2 }{ \partial t }\lambda_2 - \rho_2u_2(1-u_2)\lambda_2 + D_2( 1 - u_1 ) \nabla u_2 \cdot \nabla \lambda_2 \right) dxdt + \nonumber \\ 
&  & \int_0^T \int_{\Omega} \left( \frac{ \partial u_3 }{ \partial t }\lambda_3 + \delta_3u_3\lambda_3 - \delta_3u_2\lambda_3 + \nabla u_3 \cdot \nabla \lambda_3 \right) dxdt,\nonumber\\
&=&\left\langle \frac{\partial u }{\partial t},\lambda \right\rangle  - \left\langle F(u) , \lambda \right\rangle - \left\langle A(u) , \nabla \lambda \right\rangle, \nonumber
\end{eqnarray}
where $V = W^3$, $u,\lambda\in V$, $u = (u_1, u_2, u_3)$, $\lambda = (\lambda_1, \lambda_2, \lambda_3)$ with
\[
W = \left\{ v : v \in L^2\left(0,T;H^1_0({\Omega})\right) \mbox{ and } \frac{\partial v }{ \partial t } \in L^2(0,T; H^{-1}({\Omega}) ) \right\},
\]
and $L^2(0,T;H^1_0({\Omega})) = \left\{ v : v(x,\cdot) \in L^2((0,T)) \mbox{ and } v(\cdot,t) \in H^1_0({\Omega}) \right\}$. We use $F:V\to V^*$, $A:V \to V^*$ with
\begin{eqnarray}
\langle F(u),\lambda \rangle &=& \int_0^T \int_{\Omega} \left( u_1(1-u_1)- \delta_1u_1u_3  \right)\lambda_1 dxdt \nonumber\\
&  & + \int_0^T \int_{\Omega} \rho_2u_2(1-u_2)\lambda_2 dxdt \nonumber \\ 
&  & + \int_0^T \int_{\Omega} \delta_3 \left(  u_2- u_3 \right) \lambda_3 dxdt, \nonumber\\
\langle A(u), \nabla \lambda \rangle & = & -\int_0^T \int_{\Omega} D_2( 1 - u_1 ) \nabla u_2 \cdot \nabla \lambda_2 dxdt - \int_0^T \int_{\Omega} \nabla u_3 \cdot \nabla \lambda_3 dxdt. \label{formadebilA}
\end{eqnarray}

A weak solution $u\in V $ is a function that satisfies $\langle E(u,\delta_1),\lambda\rangle=0$ for all $\lambda \in V$.

\subsection{Formulation of the minimization problem}
As described above we propose to use an inverse problem technique in order to estimate $\delta_1$. Function $u$ represents the solution of the direct problem (the components of $u$ are the state variables of the problem) for each choice of the parameter $\delta_1$.

Let us assume that experimental information is available during the time interval $0 \leq t \leq T$. Then, the inverse mathematical problem can be formulated as:
\begin{equation}
\begin{array}{rl}
\displaystyle \mathop{ \mathrm{minimize} }_{\delta_1} & J( u , \delta_1 ) \\ \mathrm{subject \, to} & E(u,\delta_1) = 0, \\ & \delta_1 \in U_{ad},
\end{array}
\end{equation}
where the objective functional $J:V\times U_{ad}\rightarrow \mathbb{R}$ is
\begin{equation}
\label{jota1} J(u, \delta_1)=\frac{1}{2} \int_0^T \int_0^1 [u_3(x,t)-\hat{u}_3(x,t)]^2 dxdt,
\end{equation}
with $u_3(x,t)$, the excess concentration of H$^+$ ions obtained by solving the direct problem for a certain choice of $\delta_1$ and $\hat{u}_3(x,t)$, the excess concentration measured experimentally (real data). The constraints are given by $U_{ad}$, a subset of $(0,\infty)$, the set of admissible values of $\delta_1$ and $E$ is the weak formulation of the direct problem.
\begin{figure}
\begin{center}
\includegraphics[scale=0.25]{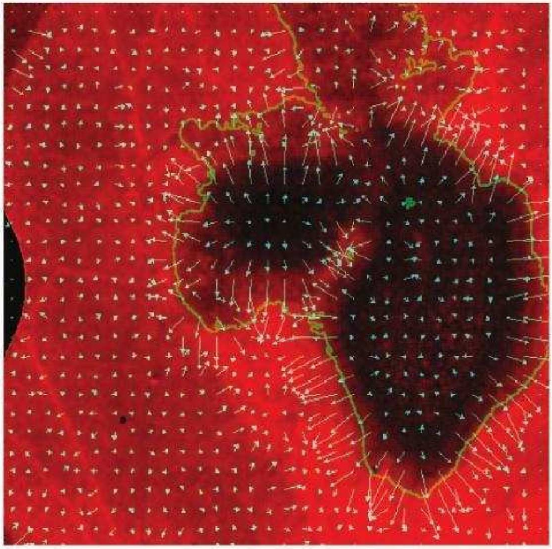}
\setlength{\belowcaptionskip}{2pt}
\caption{A map of peritumoral H$^+$ flow using vectors generated from the pH distribution around the tumor, \cite[Figure 4]{gatenby2006acid}.}
\label{fig:africa}
\end{center}
\end{figure}

We remark that, in general, there is a fundamental difference between the direct and the inverse problems. In fact, the latter is usually ill-posed in the sense of existence, uniqueness and stability of the solution. This inconvenient is often treated by using some regularization techniques
\cite{van2011source,engl1996regularization,kirsch2011introduction}.

\subsection{Formulation of the reduced and adjoint problems}
In the following, we will consider the so-called reduced problem 
\begin{equation}
\label{jotamono}
\begin{array}{rl}
\displaystyle \mathop{ \mathrm{minimize} }_{\delta_1} & \tilde{J}(\delta_1)= J(u(\delta_1),\delta_1) \\ \mathrm{subject \, to} & \delta_1 \in U_{ad},
\end{array}
\end{equation}
where $u(\delta_1)$ is given as the solution of $E(u(\delta_1),\delta_1)=0$. In order to find a minimum of the continuously differentiable function $\tilde{J}$, it will be important to compute the derivative of this reduced objective function. Hence, we will show a procedure to obtain $\tilde{J}'$ by using the adjoint approach. According to the theory exposed in \cite{brandenburg2009continuous,hinze2009optimization}, the derivative of $\tilde{J}$ is given by
\begin{equation} \label{adjoint2} \tilde{J}^{\ \prime}(\delta_1) =
\frac{\partial J}{\partial \delta_1}(u(\delta_1),\delta_1) + \left(
\frac{\partial E}{\partial \delta_1}(u(\delta_1),\delta_1)
\right)^* \lambda,
\end{equation}
where $\lambda$ solves the so-called adjoint problem
\begin{equation}\label{adjoint}
\frac{\partial J}{\partial u}(u(\delta_1),\delta_1)+\left(\frac{\partial E}{\partial u}(u(\delta_1),\delta_1)\right)^*\lambda =0.
\end{equation}

Notice that in order to obtain $\tilde{J}^{\ \prime}(\delta_1)$ we need first to compute $u(\delta_1)$ by solving the direct problem, followed by the calculation of $\lambda$ by solving the adjoint problem. For computing the second term of (\ref{adjoint2}) it is not necessary to obtain the adjoint of $\frac{\partial E}{\partial\delta_1}(u(\delta_1),\delta_1)$ but just its action over $\lambda$.

Thus, the adjoint problem (\ref{adjoint}) consists in finding $\lambda \in V $ satisfying
\begin{eqnarray}
0&=&\left\langle \frac{\partial J}{\partial u}(u(\delta_1),\delta_1),\eta \right\rangle+ \left\langle \frac{\partial E}{\partial u}(u(\delta_1),\delta_1)\eta,\lambda \right\rangle \nonumber \\
& = & \int_0^T\int_\Omega \left(- \frac{ \partial \lambda_1 }{ \partial t }\eta_1 - \eta_1( 1 -  2u_1)\lambda_1 + \delta_1\eta_1u_3\lambda_1 - D_2 \eta_1 \nabla u_2 \cdot \nabla \lambda_2\right) dxdt + \nonumber \\
& & \int_0^T\int_\Omega \left(  - \frac{ \partial \lambda_2 }{ \partial t }\eta_2 - \rho_2\eta_2( 1 - 2u_2 )\lambda_2  + D_2( 1 - u_1) \nabla \lambda_2 \cdot \nabla  \eta_2  - \delta_3\eta_2\lambda_3 \right) dxdt + \nonumber \\
& & \int_0^T\int_\Omega \left( - \frac{ \partial \lambda_3 }{ \partial t }\eta_3 + \delta_3\eta_3\lambda_3 + \nabla \lambda_3  \cdot \nabla \eta_3 + \delta_1u_1\eta_3\lambda_1 \right) dxdt + \int_0^T\int_\Omega \eta_3( u_3 - \hat{u}_3)dxdt \nonumber \\
& = & \left\langle - \frac{\partial \lambda}{ \partial t }, \eta \right\rangle + \left\langle \mathcal{H}(\lambda),\eta \right\rangle, \label{definiciondeH}
\end{eqnarray}
for all $ \eta \in V $ and $\lambda(x,T) = 0$. Then, since $\frac{\partial J}{\partial \delta_1}=0$, (\ref{adjoint2}) can be written as
\begin{equation} \label{derivadaJ}
\tilde{J}'(\delta_1)= \left(\frac{\partial E}{\partial \delta_1}(u(\delta_1),\delta_1)\right)^*\lambda = \int_0^T \int_\Omega u_1u_3\lambda_1 dxdt.
\end{equation}

\section{Designing an algorithm to solve the minimization problem}

It is worth stressing that obtaining model parameters via minimization of the objective functional $\tilde{J}$ is in general an iterative process requiring the value of the derivative. To compute $\tilde{J}'$ we just solve two weak PDEs problems per iteration: the direct and the adjoint problems. This method is much cheaper than the sensitivity approach \cite{hinze2009optimization} in which the direct problem is solved many times per iteration. We develop an implementation in MATLAB that solves the direct and adjoint problems. We use the splitting method in order to separate the direct problem in two new problems. The first one consists in a system of ordinary differential equations that contains the reaction terms of the original PDE. The second one is a PDE that contains the diffusion terms of the original PDE. The ODEs are solved by using the Runge-Kutta method ({\tt ode45} MATLAB built-in function). Since we have an ODE system for each spatial point, its resolution can be parallelized accelerating the 
time execution. The PDE is solved by the Adaptive Finite Element Method. In the next subsections we will explain the splitting method and the adaptive procedure of the FEM. It is well-known \cite{nocedal2006numerical} that gradient-based optimization algorithms require the evaluation of the gradient of the functional. The optimization problem is solved by using a Sequential Quadratic Programming (SQP) method , using the built-in function {\tt fmincon}.  

\begin{figure}
\centering
\includegraphics[scale=0.55]{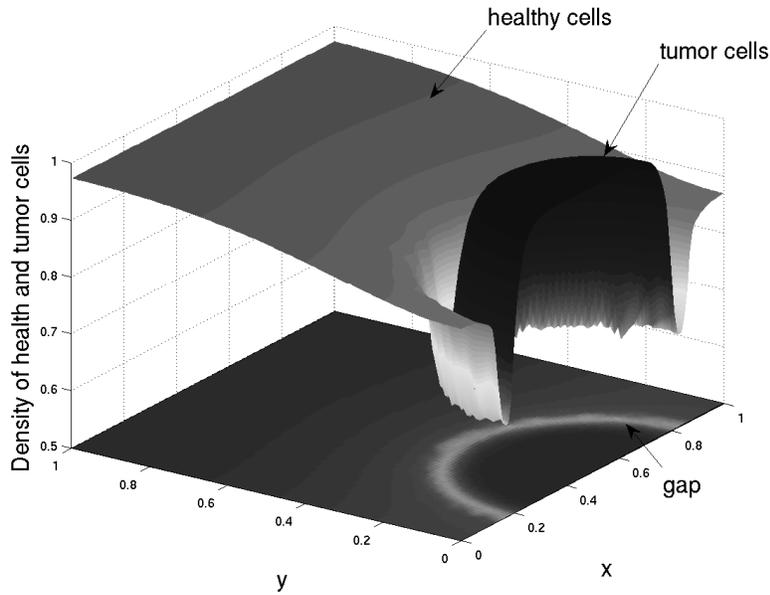}
\caption{In this figure we plot in 3D the density of the tumor and health cells that are bigger than 0.64 and its projection on, for $\delta_1 = 12.5 $ and $ t = 10$. We mark the gap, that is the region produced by the concentration of the acid that kill the health cells.}
\label{fig:directS_1}
\end{figure}

For the direct problem, Figure \ref{fig:directS_1} shows the density of health cells, tumor cells and excess concentration of H$^+$ at fixed time ($t=10$) in terms of $x$ variable.

\subsection{Solving the direct problem}
\subsubsection{Splitting method}
\begin{figure}
\begin{center}
\includegraphics[scale=0.4]{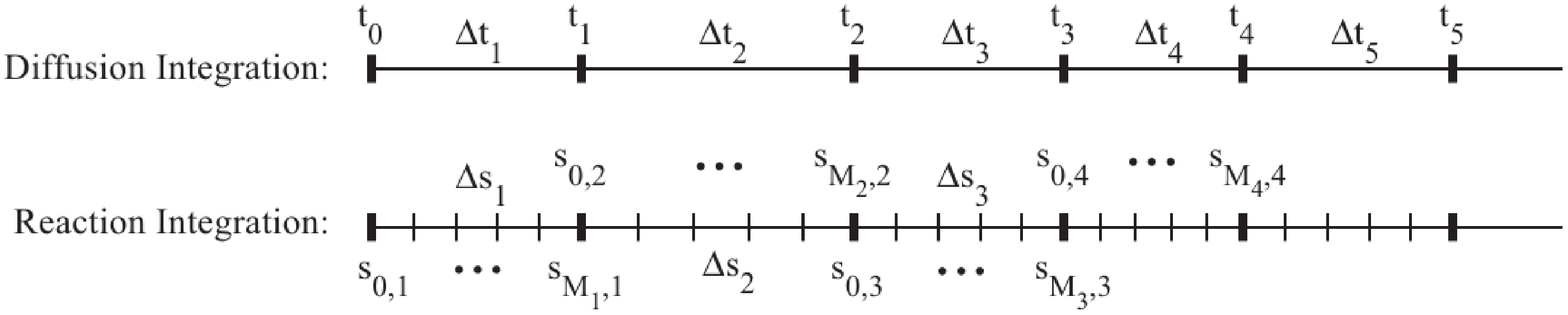}
\setlength{\belowcaptionskip}{2pt}
\caption{\cite[Section 2]{estep2008posteriori}.}
\label{fig:mos}
\end{center}
\end{figure}

\paragraph{A multiscale operator splitting}
We proceed like in \cite[Section 2]{estep2008posteriori}. For the time discretization, we introduce a theoretical framework in which each component (the reaction component $u^r$ and the diffusion component $u^d$) is solved exactly. We define a piecewise continuous approximate solution:
\[
u(x,t) = \frac{t_n - t}{\tau_n}u^{n-1}(x) +\frac{t - t_{n-1}}{\tau_n}u^{n}(x)
\]
for $t_{n-1} \leq t \leq t_n $, with the nodal values $u^n(x)$ obtained from the following procedure. We first discretize $[0, T ]$ into $0 = t_0 < t_1 < \ldots < t_N = T$ with diffusion time step $\tau$, $\tau = t_n - t_{n-1}$ for $n=1,\ldots, N$. For each diffusion step, we choose a (small) time step $\tau_{s_n} = \tau/M_n$ where $M_n \in \mathbbm{N}$, with $\tau_s = \max_{1 \leq n \leq N} \{ \tau_{s_n} \}$, and the nodes $t_{n-1} = s_{0,n} < s_{1,n} < \ldots < s_{M_n,n} = t_n$ (see figure \ref{fig:mos}). We associate the time intervals $I_n = (t_{n-1}, t_n ]$ and $I_{m,n} = (s_{m - 1,n}, s_{m,n} ]$ with these discretizations.

\subsubsection{Adaptive FEM}
\label{afem}
The adaptive procedure for FEM consists in a four step loop: (a) solve the PDE using the FEM discretization, (b) estimate \emph{a posteriori} error $\eta$ of the discrete solution, (c) mark the elements to be refined according to the relative error size of the \emph{a posteriori} error, and (d) refine the marked elements keeping the mesh conformity (for more details see \cite{nochetto2009theory}).  

Given a mesh $\mathcal{T}_n$ at time $t_n$, the element residual $R_T(u^n)$ and the jump residual $J_S(u^n)$ are defined as:
\begin{eqnarray}
R_T(u^n) & = & \frac{u^n- u^{n-1}}{\tau} -\mathcal{A}(u^n)- F(u^n), \quad T\in\mathcal{T}_n \\
J_S(u^n) & = & - A(u^{n+}) \cdot \nu^+ - A(u^{n-}) \cdot \nu^-, \quad S\in\mathcal{S}_n 
\end{eqnarray}
were $\mathcal{S}_n $ are the edges of $\mathcal{T}_n$ and $\mathcal{A}$ is the strong form of the $A$ operator as defined in (\ref{formadebilA}). 

We define the local error indicator $\eta(T)$ by
\[
\eta(T)^2 = 
H_T^2 \|R_T(u^n)\|_{L^2(T)}^2 + \sum_{S\in \partial T} H_S \|J_S(u^n)\|_{L^2(S)}^2,
\]
were $H_T$ is the diameter of $T$ and $H_S$ is the length of the edge $S$. If $S$ is an edge of an element, then
\[
\eta(S)^2=H_S \|J_S(u^n)\|_{L^2(S)}^2.
\]

The residual-type error estimator of $\Omega$ with respect to the mesh $\mathcal{T}_n$ is
\[
\eta(\Omega)^2 = \sum_{T\in\Omega}\eta(T)^2. 
\]

\subsubsection{Algorithm}\label{algdirect}
\begin{itemize}
\item[] STEP 0: Set an initial condition $u^0(x) = u(x,0)$ on the coarse uniform mesh $\mathcal{T}_0$. Set $\varepsilon_{TOL} > 0$.
\item[] STEP 1: Given $u^{n-1}(x)$ do the following steps to compute $u^n(x)$ if $n \leq N$.
\item[] STEP 2: Compute $u^r(x,t)$ satisfying the reaction equation:
\begin{eqnarray}
\left\langle \frac{\partial u^r}{\partial t} , \phi \right\rangle_{I_n} 
&=& \left\langle F(u^r),\phi \right\rangle_{I_n}  \nonumber\\
u^r(x,t_{n-1}^+) &=& u^{n-1}(x) \nonumber	
\end{eqnarray}
for $s_{0,n} < t \leq s_{M_n,n} $  and for all $\phi\in V$.
\item[] STEP 3: Compute $u^d(x,t)$ satisfying the diffusion equation:
\begin{eqnarray}
\left\langle \frac{\partial u^d}{\partial t} , \phi \right\rangle_{I_n} 
& = &  \left\langle A(u^d) , \phi \right\rangle_{I_n} \nonumber\\
u^d(x,t_{n-1}^+) &=& u^r(x,t_n) \nonumber
\end{eqnarray}
for $t_{n-1} < t \leq t_n $ and for all $\phi\in V$. Set $u^n(x) = u^d(x,t_n)$.
\item[] STEP 4: Compute the \emph{a posteriori} error $\eta(\Omega)$. If $\eta(\Omega) < \varepsilon_{TOL}$, set $n=n+1$ and go to STEP 1.
\item[] STEP 5: Mark and refine, and go to STEP 2.
\end{itemize}

In STEP 2, we compute $u^r$, the reaction component, by using the {\tt ode45} MATLAB  built-in function for each node of the current mesh, allowing a parallelization strategy. In STEP 3, the diffusion component is solved by using FEM. In STEP 5, according to \cite{ffw}, we use the {\tt bulk} algorithm to mark and the {\tt RedGreenBlue} algorithm  to refine. The bulk algorithm  defines the set $\mathcal{E}$ of marked edges such that
\[
\sum_{E\in\mathcal{E}} \eta(E)^2 \geq \theta \sum_{S\in\mathcal{S}_n} \eta(S)^2, 
\]
or it contains all the edges of marked elements $T\in \mathcal{K} \subset \mathcal{T}_n$ that satisfy
\[
\sum_{K\in\mathcal{K}} \eta(K)^2 \geq \theta \sum_{T\in\mathcal{T}_n} \eta(T)^2. 
\]
where $\mathcal{S}_n$ is the set of edges of $\mathcal{T}_n$ and $\theta \in [0,1]$.

\subsection{Solving the adjoint problem}
In order to solve the adjoint problem we shall use FEM. The spatial discretization is the coarse mesh used for the initial mesh in the direct problem. We denote $\lambda^n(x)=\lambda(x,t_n)$ for $n=0,\ldots,N$.

\subsubsection{Algorithm}\label{algadj}
\begin{itemize}
\item[] STEP 0: Set the final condition $\lambda^N(x)=\lambda(x,T) = 0$ on the initial mesh $\mathcal{T}_0$.
\item[] STEP 1: Given $\lambda^{n}(x)$ do the following steps to compute $\lambda^{n-1}(x)$ if $n \geq 1$.
\item[] STEP 2: Do an implicit Euler step in time, and FEM in space to approximate the adjoint variable $\lambda^{n-1}$ by solving the linear system $\lambda^{n-1} - \lambda^n - \tau K(\lambda^{n-1})=0$, where $K$ is the discretization of $\mathcal{H}$ as defined in (\ref{definiciondeH}).
\item[] STEP 3: Set $n=n-1$ and go to STEP 1.
\end{itemize}

\subsection{Solving the minimization problem}
The \texttt{fmincon} MATLAB built-in function was used to solve the minimization problem. The chosen algorithm in the \texttt{fmincon} function was the trust-region-reflective method, where the derivative of the objective function $\tilde{J}$ was computed according to \ref{derivadaJ}.

\subsubsection{Algorithm}\label{algmin}
The method we will use for minimizing the functional $\tilde{J}$ can be summarized as follows:
\begin{itemize}
\item[] STEP 0: Give an initial guess $\delta_1^0$ for the parameter.
\item[] STEP 1: Call the {\tt fmincon} function and obtain the solution $\delta_1^*$, providing the objective function $\tilde{J}(\delta)$ and its derivative $\tilde{J}'(\delta)$ according to (\ref{jotamono}) and (\ref{derivadaJ}), respectively. 
\end{itemize}

In order to compute $\tilde{J}(\delta)$ and $\tilde{J}'(\delta)$ is necessary to solve the direct and adjoint problems.

\section{Numerical experiments}
The goal of this section is to test and evaluate the performance of an adjoint-based optimization method, by executing some numerical simulations of Algorithm \ref{algmin} for some test cases. The experiments were run in MATLAB, in a PC running Linux, Intel(R) Core(TM) i7-3770K CPU, 3.50GHz. 
 
First consider an optimization problem that consists in minimizing the functional defined in (\ref{jotamono}), where $\hat{u}_3(x,t)$ is generated via the direct problem for some $\hat{\delta}_1$ with the choice of model parameters $\rho_2 = 1$, $D_2= 4\times 10^{-5}$ and $\delta_3 = 1$. We choose several values of $\hat{\delta}_1$, for instance $\hat{\delta}_1 = 0.5, \, 4,\, 12.5 ,\, 16 $, because they show a different behavior of tumor invasion, according to \cite{gatenby1996reaction}.

\begin{figure}[!ht]
\centering
\includegraphics[scale=0.42]{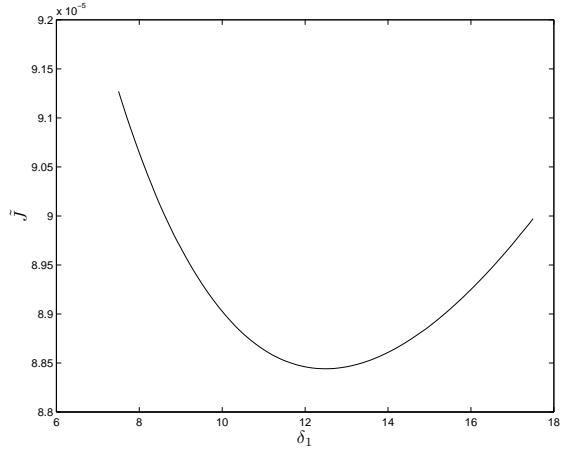}
\caption{The functional $\tilde{J}$ for $\hat{u}_3$ generated with $\hat{\delta}_1=12.5$.}
\label{fig:figurej}
\end{figure}

Figure \ref{fig:figurej} shows the value that the functional defined in (\ref{jotamono}) takes for different values of $\delta_1$, for $\hat{u}_3$ generated with $\hat{\delta}_1=12.5$. It is worth mentioning that, even when we do not know in advance if the optimization problem has a unique solution, $\tilde{J}$ looks convex with respect to $\delta_1$.

\begin{figure}
\begin{tabular}{cc}
\includegraphics[scale=0.45]{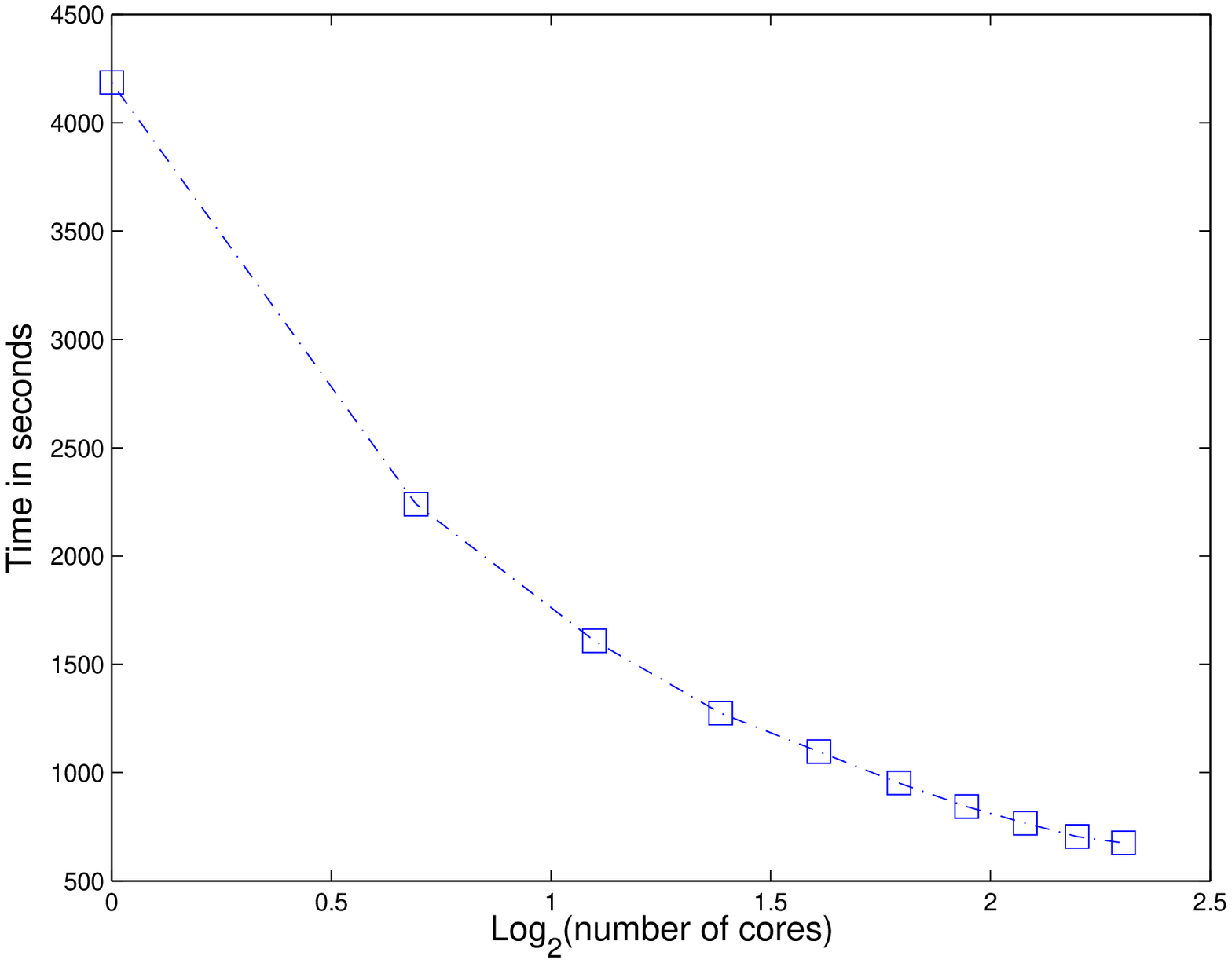} & \includegraphics[scale=0.45]{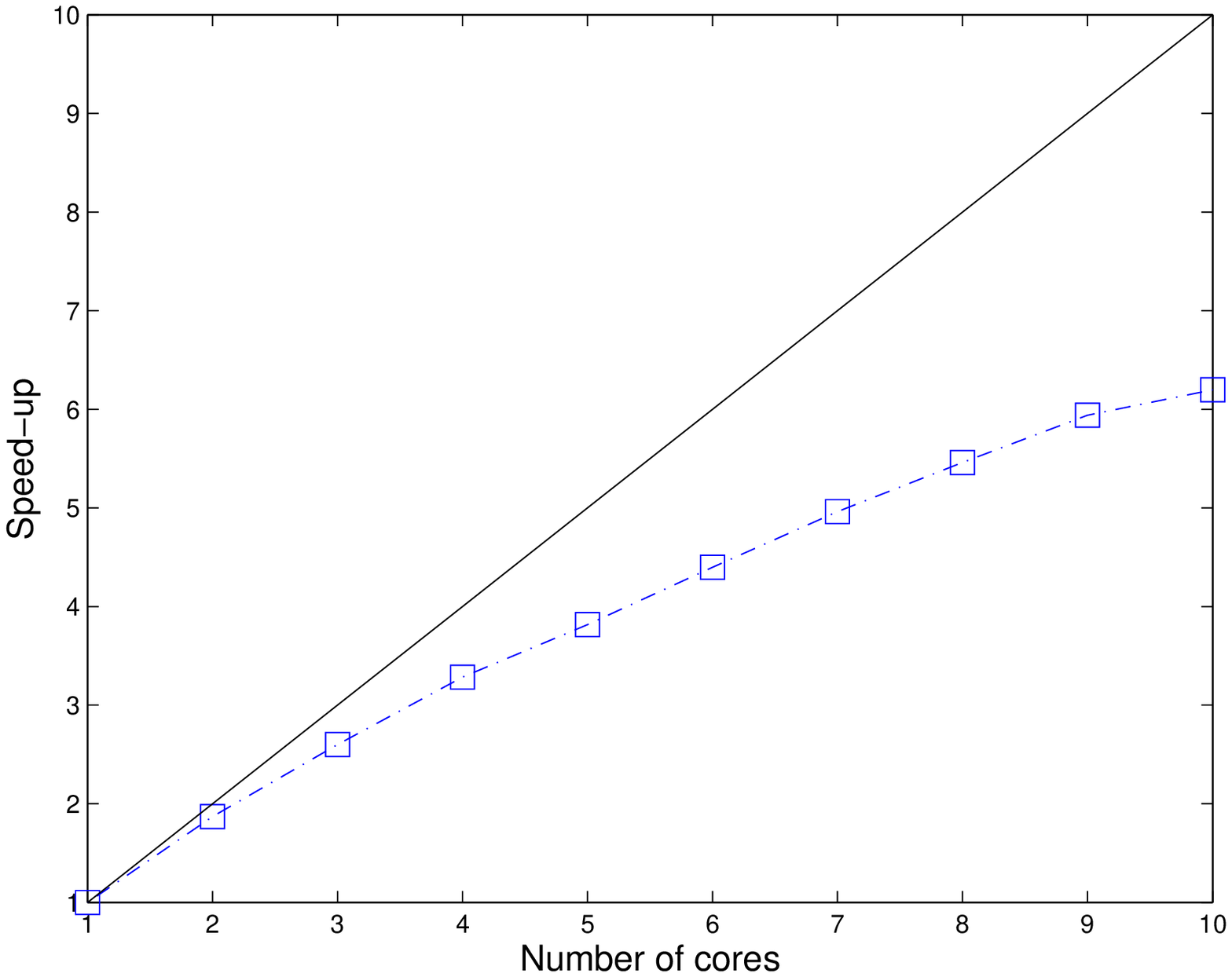} \\
(a) & (b)
\end{tabular}
\caption{Time in seconds (a), and speed-up (b) of the execution of the direct problem.}
\label{fig:figurecoretime}
\end{figure}

The idea of these test cases is to investigate how close the original value of the parameter can be retrieved (even in the presence of noise), and how efficiently these computations can be done. Regarding computational efficiency, one of the most expensive parts of Algorithm \ref{algmin} is the resolution of a system of ODEs (STEP 2 of Algorithm \ref{algdirect}). Since we have a system of ODEs for each node of a current mesh, a parallel strategy (each processor solves a system of ODES for one node) is the best and natural option to reduce time execution. For this particular inverse problem, parallelization is not an option, but it is a need, because the direct problem could be called many times by the optimization solver. For example, figure \ref{fig:figurecoretime}(a) shows how many seconds takes to solve the direct problem. Figure \ref{fig:figurecoretime}(b) shows the speed-up. In addition, since we have to go through all the nodes in order to compute the \emph{a posteriori} error, we have also parallelized this computation.

On the other hand, the method used for the minimization algorithm (\ref{algmin}) is the trust-region-reflective method \cite{coleman1996reflective,matlab}, where the option {\tt GradObj} is on by default (the gradient of the objective function must be supplied). If we would use another algorithm (like the active-set or the interior-point algorithms) where the gradient can be estimated by finite differences, then the computational cost of this algorithm makes it not practical. That is the reason for which we compute the exact derivative of the functional using the adjoint method.

We have run the Algorithm \ref{algmin} for different values of $\hat{\delta}_1$ taking the initial condition $\delta_1^0$ randomly. Averaging the different solutions and taking the standard deviation of all of these experiments, the retrieved parameter is obtained very accurately.

The algorithmic parameters for Algorithm \ref{algdirect} and \ref{algadj} are: $\tau = 0.1$, $T = 10$, the initial coarse mesh $\mathcal{T}_0$ has $512$ triangular elements, $\varepsilon_{TOL}=10^{-5}$, and $\theta = 1/2$.

The algorithmic parameters for Algorithm \ref{algmin} are: the feasible set for the optimization problem (\ref{jotamono}) is $U_{ad}=[0,20]$, the method used is trust-region-reflective, the option {\tt GradObj} is on, and the maximum of function evaluations is 100.

\begin{table}
\centering
\begin{tabular}{ r | c | c }
$\hat{\delta}_1$ & $\bar\delta_1$ & $ S $ \\ \hline
4    & 4.2666  & $\pm$ 7.0640 $\times10^{-3}$ \\
12.5 & 12.4937 & $\pm$ 7.1875 $\times10^{-4}$ \\
16   & 16.6246 & $\pm$ 1.8826 $\times10^{-5}$
\end{tabular}
\caption{Experiments for randomly initial data $\delta_1^0$}
\label{tabla1}
\end{table}

We emphasize that we have retrieved accurately the value of $\hat{\delta}_1$ independently of the value of $\delta_1^0$. Thus, in the next experiment we will consider a fixed value $\delta_1^0=8$.

It is well-known that the presence of noise in the data may imply the appearance of strong numerical instabilities in the solution of an inverse problem \cite{bertero2006inverse}.

One of the experimental methods to obtain values of $\hat{u}_3$ is by using fluorescence ratio imaging microscopy \cite{martin1994noninvasive}. 

As it is well-known, measurements are often affected by perturbations, usually random ones. Then we perform numerical experiments where $\hat{u}_3$ is perturbed by using Gaussian random noise with zero mean and standard deviation $\sigma = 0.1, \, 0.15, \, 0.2$.  In Tables \ref{delta1=4}-\ref{delta1=16}, for each value of $\sigma$, we show the average $\bar\delta_1$ for 10 values of $\delta^*_1$, the standard deviation $S$ and the relative error $e_{\delta^*_1}$.

\begin{table}
\centering
\begin{tabular}{ c | c | c | c }
$\sigma$ &  $\bar\delta_1$ & $ S $  & $e_{\delta_1}$  \\ \hline
0.1000 & 3.1125  & $\pm$ 0.8624   & 0.2219  \\
0.1500 & 3.5409   & $\pm$ 1.8611   & 0.1148  \\ 
0.2000 & 3.4471   & $\pm$ 2.3701   & 0.1382  \\ 
\end{tabular}
\caption{Experiments for $\hat{\delta}_1 = 4$}
\label{delta1=4}
\end{table}

\begin{table}
\centering
\begin{tabular}{ c | c | c | c }
$\sigma$ &  $\bar\delta_1$ & $ S $  & $e_{\delta_1}$  \\ \hline 
0.1000 & 11.6235  & $\pm$ 2.6314  &  0.0701 \\
0.1500 & 12.3825  & $\pm$ 4.6561  &  0.0094 \\
0.2000 & 11.9537  & $\pm$ 5.5749  &  0.0437 \\
\end{tabular}
\caption{Experiments for $\hat{\delta}_1 = 12.5 $}
\label{delta1=12.5}
\end{table}

\begin{table}[!h]
\centering
\begin{tabular}{ c | c | c | c }
$\sigma$ &  $\bar\delta_1$ & $ S $  & $e_{\delta_1}$  \\ \hline
0.1000 & 16.6996  & $\pm$  2.1280   &  0.0437 \\
0.1500 & 17.0926  & $\pm$  5.1026   &  0.0683 \\
0.2000 & 17.6308  & $\pm$  2.4753   &  0.1185 \\ 
\end{tabular}
\caption{Experiments for $\hat{\delta}_1 = 16$}
\label{delta1=16}
\end{table}

\section{Final conclusions and future work}
In this paper we have solved a parameter estimation problem following the model proposed by \cite{gatenby1996reaction} in a two-dimensional space. The inverse problem is formulated as an optimization problem in order to find the parameter $\delta_1$ (the destructive influence of H$^+$ ions in the healthy tissue). 

The direct problem was solved by the splitting technique together with Adaptive Finite Element Method, for the purpose of controlling the numerical error and defining a parallel strategy. A gradient-based method was used to solve the optimization problem. The derivative of the objective functional was computed using the solution of the adjoint problem. 

The experiments were run in MATLAB recovering several values of the parameter $\delta_1$ representing different scenarios. Also, a stability analysis was performed using random noise to simulate perturbations in the experimental data.

We consider that the results are accurately enough. In most cases the parameters are retrieved with a relative error less than $20\%$.

As a future work we propose to consider the possibility to find optimal parameters related to therapeutic procedures like in \cite{mcgillen2012applications,martin2012predicting}. 

\section*{Acknowledgments}
We appreciate the courtesy of Sebastian Pauletti, from IMAL, Santa Fe, Argentina, who strongly contributed with information above splitting method.

The work of  the authors  was partially supported by grants from CONICET, SECYT-UNC and PICT-FONCYT.

\bibliographystyle{model1-num-names}

\end{document}